\documentclass[a4paper,10pt,leqno]{article}

\usepackage{amssymb}

\usepackage{amsmath}
\usepackage{caption}
\usepackage{booktabs}
\usepackage{amsthm}
\usepackage{url}
\usepackage{xcolor}
\usepackage{algorithm}
\usepackage{algorithmicx}
\usepackage{graphicx}
\usepackage{mathtools}
\usepackage{algpseudocode} 
\usepackage{multirow}
\usepackage{array}
\usepackage{comment}
\usepackage[colorlinks, citecolor=blue]{hyperref}

\newtheorem{definition}{Definition}

\DeclareRobustCommand{\rchi}{{\mathpalette\irchi\relax}}
\newcommand{\irchi}[2]{\raisebox{\depth}{$#1\chi$}}
\newcommand{\numberset}{\mathbb}
\newcommand{\N}{\numberset{N}}
\newcommand{\Z}{\numberset{Z}}
\newcommand{\R}{\numberset{R}}

\newcommand{\I}{\numberset{I}}

\newcommand{\T}{\numberset{T}}

\DeclareMathOperator{\Id}{Id}

\title{Symmetric solutions of the $n$-body problem: a numerical study of Floquet multipliers and Morse indices}
\author{
    Diego Berti, Gian Marco Canneori, Roberto Ciccarelli \\
	Irene De Blasi, Margaux Introna, Davide Polimeni \\ [1em]
}

\date{}

\begin{document}
	\maketitle

	\begin{abstract}
		In this paper, we consider periodic solutions of the $n$-body problem that satisfy symmetry constraints, expressed through invariance under finite group actions. We focus on their stability properties and present algorithms specifically designed for the computation of Floquet multipliers and Morse indices. Numerical results are provided to illustrate our methods in both two and three dimensional configuration spaces, and for different choices on the number of bodies. 
	\end{abstract}

	\section{Introduction}
\label{sec:intro}

\noindent The $n$-body problem stands as one of the most important and challenging models in Celestial Mechanics, having attracted the attention of mathematicians and physicists for centuries. Its formulation is both clear and simple: consider $n$ massive particles in the space, subject to their mutual gravitational attraction. Starting from given initial conditions, the question is whether one can understand and predict the qualitative behaviour of their motion. Gravitational Newton laws provide the motion equations of $x_1,\ldots, x_n\in\mathbb{R}^3$ heavy bodies with masses $m_1,\ldots, m_n>0$, which read
\[
    m_k\ddot{x}_k(t) = -\sum\limits_{j\neq k} m_k m_j \frac{x_k(t)-x_j(t)}{\|x_k(t)-x_j(t)\|^3},\quad\text{for each}\ k\in\{1,\ldots,n\},
\]
where $\|\cdot\|$ is the Euclidean norm of $\mathbb{R}^3$. As a first observation, any collision between two or more bodies is a singularity of the motion equations and breaks the completeness of the associated flow.

As a matter of fact, what appears to be an elementary formulation soon reveals a remarkable complexity as the number of bodies increases. When only 2 bodies are involved, we deal with the famous \emph{Kepler problem}, which is an example of completely integrable system (i.e., starting from any initial configuration, closed formul\ae\ for solutions can be found). The regular property of integrability is broken by introducing one or more additional bodies in the picture, as remarkably noticed in the pioneering studies of Henri Poincaré \cite{poincare1893methodes}. Small changes in the initial conditions of known solutions can lead to unexpected and unpredictable phenomena, raising a reasonable conjecture on the deterministic chaos of the $n$-body problem (see \cite{Moser, Bolotin2006, LlibreSimo, GuardiaMartinSeara,Boekholt, Leigh}). Moreover, interesting open questions on the long term behaviour and qualitative description of such systems still persist (see the book \cite{Montgomery2024}). 

Henri Poincaré conjectured that periodic orbits are dense in the 3-body problem and claimed that they play a central role in capturing the system's complexity. In this picture, periodic orbits provide a systematic and useful tool to approximate increasingly intricate trajectories, serving as empirical validation of topological chaos. For this reason, the community has shown great interest in finding periodic solutions for the $n$-body problem, employing both numerical and analytical approaches. In particular, since the 90s, variational methods have been successfully applied to produce collision-free periodic solutions (see \cite{Bahri, Serra_Terr_92, MajerTerracini1, MajerTerracini2, MajerTerracini3}). In short, one can consider the Lagrange action functional associated with the motion equations 
\[
    \mathcal{A}(x) = \int_0^T\left[\frac12\sum\limits_{i=1}^n m_i\|\dot{x}_i(t)\|^2+ \sum\limits_{i<j}\frac{m_i m_j}{\|x_i(t)-x_j(t
    )\|}\right]\,dt
\]
and look for its critical points, which naturally correspond to solutions of the $n$-body problem through a least action principle. Among all periodic solutions, one can focus on those who satisfy a symmetry constraint, especially if this restriction makes the search of critical points more feasible. The first successful effort for the planar $n$-body problem traces back to \cite{Bessi}, where the authors proved the existence of non-collision solutions using variational methods. Later on, in the celebrated paper \cite{chen_mont00}, Alain Chenciner and Richard Montgomery proved the existence of the spectacular figure-eight solution, previously discovered in \cite{Moore}. Subsequently, Davide Ferrario and Susanna Terracini enlarged the range of applicability of variational methods to symmetric $n$-body problems in the celebrated paper \cite{FerrarioTerracini}. Discovering elementary algebraic conditions, they allow to consider a wide class of symmetries, in which periodic solutions can be found. The symmetry constraint on the configurations of the bodies is expressed by using the action of a finite group $G$, and they look for critical points of the action functional among $G$-equivariants loops. Notably, variational methods produced many other results in the last two decades, enlarging the set of symmetric periodic solutions for the $n$-body problem (see, e.g., \cite{Chen,Ramos, ChencinerVenturelli, Barutello_2004, gronchi, Sim}).)

Back to the paper \cite{FerrarioTerracini}, alongside with a huge theoretical effort to deal with collisions between the bodies and to provide a solid theory on $G$-equivariant loops, the authors conceived a useful algorithm to produce numerical symmetric solutions. Their theoretical results found a practical confirmation and the software \textit{symorb} was created (\cite{symorb, Ferrario2024}). It consists of a numerical pipeline based on a combination of Python, Fortran and GAP, which allows to choose a finite group, consider the symmetry involved and look for equivariant critical points of the action functional. 

Recently, a new version of the original software was presented in \cite{SymOrbJL_github, symorbJL}. Distributed as a Julia package, \emph{SymOrb.jl} represents a profound and modular redesign of the previous implementation. This new architecture provides a flexible and extensible framework, enabling seamless integration of both quantitative and qualitative tools, such as stability analysis and topological index computations.
It also introduces more efficient optimization routines, improving computational performance.
Finally, the new structure of \emph{SymOrb.jl} enables the systematic organization of symmetric orbits into databases and supports advanced numerical methods.

Starting from a collection of numerical results produced with \textit{SymOrb.jl}, the present paper aims to establish a set of benchmark cases, by examining their stability properties under variations of the action and small perturbations of the orbits' initial conditions. To this end, two classical stability indicators are considered: the Morse index and the Floquet multipliers (\cite{Morse1934, Smale1965, teschl2024ordinary}). We also outline the main steps required to compute a numerical approximation of these indicators in our setting.

Overall, this work offers a first comprehensive illustration of how \textit{SymOrb.jl}, complemented by its stability modules, can be employed to find highly accurate numerical approximations of solutions to the $n$-body problem and to analyze their stability through different methodological approaches.

\paragraph{Mathematical setting and problem statement} 
For $n\geq2$ particles, we denote their masses as $m_1,\dots,m_n>0$ and their positions as $x_i\in \R^d$, where $d=2,3$. Since the $n$-body problem is a mechanical system, we can introduce the potential function 
\[
    U(x_1, \dots, x_n)=\sum_{j<i}\dfrac{m_jm_i}{\|x_j-x_i\|},\quad\|\cdot\|\ \text{Euclidean norm of}\ \R^d,
\]
so that the equation of motion of the $i-$th particle is given by 
\begin{equation}\label{eq:nbody}
    m_i \ddot{x_i}(t)=\dfrac{\partial U}{\partial x_i}(x_1(t), \dots, x_n(t)).
\end{equation}
It is easy to see that the centre of mass of the system is invariant under translations, so it is useful to fix it at the origin and to introduce the configuration space 
\[
\rchi \coloneq \left\{x=(x_1,\dots,x_n)\in(\R^d)^n    \ :\  \sum_{i=1}^n m_i x_i = 0\right\}.
\]
Note that the potential function $U$ has a singularity whenever $x_i=x_j$ for some $j\neq i$. The singularity set of $U$ can be described as follows
\[
\Delta_{ij}=\left\{x\in \rchi \ :\ x_i=x_j\right\}, \quad \Delta=\bigcup_{i,j}\Delta_{ij}, 
\]
and, as a consequence, the set of non-colliding configurations is nothing but $\hat \rchi=\rchi\setminus \Delta$. We also define the kinetic energy as 
\[
K(\dot x)={\frac{1}{2}}\sum_{i=1}^nm_i\|\dot x_i\|, \quad \dot x=(\dot x_1, \dots, \dot x_n)\in T_x\rchi, 
\]
where $T_x\rchi$ denotes the tangent space to $\rchi$ at $x\in\rchi$. In this way, the Lagrangian function associated to equations \eqref{eq:nbody} can be expressed as 
\[
L(x, \dot x)=U(x) + K(\dot x).
\]
Classical periodic solutions of \eqref{eq:nbody} are trajectories $x(t)=(x_1(t), \dots, x_n(t))\in \hat\rchi$ such that $x(T)=x(0)$ for a suitable $T>0$: the smallest $T$ that satisfies this condition is called the \emph{period} of the orbit. Given $T>0$, let us then define the torus of length $T$ as $\T=\R/T\Z$, and consider the space of $H^1$ \emph{T-periodic loops} (possibly with collisions between the bodies) 
\[
\Lambda = H^1(\mathbb{T};\rchi)=\left\{x, \dot x\in L^2\left([0, T], \rchi\right) \ :\ x(0)=x(T)\right\}.
\]
We then denote by
\[
    \hat{\Lambda} = H^1(\mathbb{T};\hat{\rchi})\subset \Lambda 
\]
the open sub-set of collision-less loops. The Lagrange action functional on $\Lambda$ reads
\begin{equation}
\label{action functional}
    \begin{aligned}
    \mathcal{A}(x) &= \int_0^T L(x(t),\dot{x}(t))\,dt  \\
    &= \int_0^T\left[\frac12\sum\limits_{i=1}^n m_i\|\dot{x}_i(t)\|^2 + \sum\limits_{i<j}\frac{m_i m_j}{\|x_i(t) - x_j(t)\|}\right]\,dt
    \end{aligned}
\end{equation}
and has the following property: if $x\in\hat{\Lambda}$ is a critical point of $\mathcal{A}$, then $x$ is a $T$-periodic solution of \eqref{eq:nbody} (see \cite{FerrarioTerracini}). This is the variational principle we are going to refer to in order to find $T-$periodic solutions of our system.

\paragraph{Group actions, G-equivariance and optimisation} 

\noindent In this paragraph we briefly recall the basic definitions needed to  impose symmetry constraints in the bodies motions. Our main reference here is \cite{FerrarioTerracini}, while the interested reader can find a detailed account on G-equivariance and group actions in \cite{symorbJL}. 

Let $G$ be a finite group; we say that $G$ acts on a set $X$ if there exists a map $\phi:G\times X\to X$ such that 
\[
\phi(1,x)=x\quad \forall x\in X, \quad \phi(g,\phi(h,x))=\phi(gh, x) \quad \forall g,h\in G, \  x\in X. 
\]
In the context of $T-$periodic loops in $\rchi$, we can describe the action of a finite group $G$ by defining how its elements behave on the space, time and body labels. More precisely, given $x\in \Lambda$ and $t \in \T$, it is possible to represent the action of an element $g\in G$ over $x(t)$ through the following homomorphisms
\[
\rho: G\to O(d), \quad \tau: G\to O(2) , \quad \sigma: G \to \Sigma_n,  
\]
where $O(d)$ and $O(2)$ denote the orthogonal groups of dimension $d$ and $2$ respectively, and $\Sigma_n$ is the group of all permutations of $n$ objects labelled in $\{1, \dots, n\}$.

In practice, given $g\in G$, $\rho(g)$ describes how $g$ acts on the $d-$dimensional space where every $x_i$ lies, while $\sigma(g)$ describes the possible interchanging of bodies along the loop. As for $\tau(g)$, it is used to represent possible symmetries or recurrences of the orbit over a period $T$. In Section \ref{sec:numeric}, a group $G$ will be identified by its generators, whose representations on $O(d), O(2)$ and $\Sigma_n$ will be provided. 

Given a finite group $G$, we can define the set of the (possibly colliding) \emph{G-equivariant} loops in $\Lambda$ as the sets
\[
\Lambda^G=\left\{x\in \Lambda \ : \ (gx)(t)=x(t), \ \forall t\in \T, \ g\in G\right\}, \quad \hat\Lambda^G=\left\{x\in \Lambda^G \ :\ x\in \hat\Lambda\right\}
\]
that is, the set of all elements of $\Lambda$ (resp. $\hat\Lambda$) which are invariant under the action of $G$. It is possible to show (see \cite[Lemma 2.4]{symorbJL}) that, if $\bar x\in \hat \Lambda^G$ is a critical point of the restriction $\mathcal A_{|_{\Lambda^G}}$, then it is a critical point of $\mathcal A$ over the whole $\Lambda$, and hence a periodic $G$-equivariant solution of \eqref{eq:nbody}. With suitable constraints on $G$, one can ensure \emph{a priori} the existence of such collision-less critical points (see again \cite{FerrarioTerracini, symorbJL} for more details). 

To find critical points of $\mathcal A_{|_{\Lambda^G}}$, it is possible to operate a further reduction by considering a particular subinterval $\mathbb I\subset \T$, called the \textit{fundamental domain}. In short, the fundamental domain $\mathbb{I}$ is such that 
\begin{enumerate}
    \item[(A)] defined $\bar G=G/\ker\tau$, one has 
    \[
    \T = \bigcup_{[g]\in \bar G}\tau(g^{-1})\mathbb I, \quad |\mathbb I|=\dfrac{|\T|}{|\bar G|}
    \]
    \item[(B)] if $\bar y: \mathbb I \to \rchi$ is a critical point of the restricted action 
    \[\mathcal A_{\mathbb I}(y)=\int_{\mathbb I}L(y(t), \dot y(t))dt\]
    over a suitable set of fixed-ends trajectories, then the symmetrised path given by the concatenation of $g \bar y,$ $g\in \bar G$,
    is a solution of the $G-$equivariant $n-$body problem (see \cite[Theorem 3.2]{symorbJL}). 
\end{enumerate}
Some considerations on properties (A) and (B) are in order to understand the optimisation process and, later on, the computation of the stability indicators. Property (A) claims that it is possible to re-construct the whole period simply looking at the action of $\tau(G)$ on the segment $\mathbb I$.  As for property (B), it allows to restrict our search for critical points of the action in the set of segments of our $T-$periodic orbits with suitable constraints at the endpoints. In particular, $y$ must belong to a class $Y$ of fixed-ends-type paths  whose precise definition depends on the form of $\tau(G)$ (see \cite[Proposition 3.1]{symorbJL}).
 Without loss of generality, we will take $\I = [0,\pi]$ and $T = l\pi$, $l\in\N$.

The optimisation of the restricted action functional $\mathcal A_{\I}$  is the key problem tackled by the numerical algorithm, which is the core of \textit{SymOrb.jl}. Given a symmetry group, and, consequently, the corresponding fundamental domain $\I$ and the set of the paths $Y$, the routine takes advantage of a wide variety of optimisation and refinement methods to compute an approximation of a critical point $\bar y$ of $\mathcal A_{\I}$ over $Y$. Finally, by symmetrising $\bar y$, an approximated solution $\bar x$ of \eqref{eq:nbody} over the period $T$ is obtained. The output of the optimisation is provided in two ways: first, the Fourier coefficients of $\bar y$ on $\I$, and secondly, the pointwise representation of $\bar x$ on $\T$. 

Once the approximated solution is provided, one can proceed with its stability analysis, which is the topic of next Section.

\section{Stability indicators} \label{sec:stab}

Once an approximated solution of \eqref{eq:nbody} is obtained via the \textit{SymOrb.jl} algorithm, it is possible to study its stability as a periodic orbit of the $n-$body problem. In this section, we will propose two different ways to evaluate stability. First of all, we will define the Floquet multipliers, that can be used to evaluate the change rate of a periodic orbit under small changes in its initial condition. Secondly, we will introduce the Morse index, which is related to the variation of the Lagrangian action value when the orbit is slightly deformed. For both quantities, we will provide the basic definitions and the computing methods in our context; the interested can find more details on such stability indicators in \cite[Chapters 3 and 12]{teschl2024ordinary}, \cite{Morse1934, Smale1965}, as well as other applications of index theory to problems in Celestial Mechanics in \cite{AmbrosettiCotiZelati1993, Rabinowitz1978, BarutelloHu}. 

\subsection{Floquet multipliers}

 Let us start by taking a general ordinary nonlinear differential equation of the form 
 \begin{equation}\label{eq:init prob}
     \dot x(t) =f(x(t))
 \end{equation}
where $x: \mathbb R \to \mathbb R^N$, for some $ N\ge 1$ and $f: \mathbb R^N \to \mathbb R^N$. Suppose that the above system admits a $T-$periodic solution, $T>0$, denoted by $\bar x(t)$.  The linear stability of such periodic solution can be analysed by passing to the linearised system 
\begin{equation}\label{eq:linearized}
    \dot y (t)=A_{\bar x}(t)y(t), \ A_{\bar x}(t)=D(f(\bar x(t)):   
\end{equation}
clearly, $A_{\bar x}(t)$ is a periodic $N\  \times \  N$ matrix and depends on the particular solution $\bar x(t)$. 
Following the theory of periodic linear systems, it is possible to define the \emph{principal matrix solution} $X$ as the unique solution of the differential system
\begin{equation}
\label{eq:matrix_eq}
    \begin{cases}
    \dot X (t) = A_{\bar x}(t)\,X(t)
    \\
    X(0) = I_N,
\end{cases}
\end{equation}
where $X(t) \in \mathbb R^{N \times N}$ for every $t\in \mathbb R$ and $I_N$ is the $N-$dimensional identity. Starting from $X(t)$ one can define the \emph{Monodromy matrix} related to problem \eqref{eq:linearized} as $
{\mathcal M}_{\bar x} = X(T)$.

The spectral properties of the monodromy matrix are related to the linear stability of the periodic solution $\bar x$: if all eigenvalues of ${\mathcal M}_{\bar x}$ have complex modulus less or equal to $1$, then $\bar x$ is \emph{linearly stable}. Otherwise, $\bar x$ is  \emph{linearly unstable}. The eigenvalues of $M_{\bar x}$ are called  the \emph{Floquet multipliers} associated to the orbit $\bar x$ (see, for instance, \cite[Theorem 12.4]{teschl2024ordinary} \cite{Barut16} and  \cite{calleja}). 

It is important to note that the analytical computation of Floquet multipliers is, in general, a challenging task, with only a limited number of known results. As a result, numerical methods are typically preferred for their determination (see the next sections for references in the case of $n$-body problem).

\subsubsection{Computation of Floquet multipliers for the $n$-body problem}

Following previous works in the literature (see, e.g., \cite{neworbits, Barut16, calleja}), we analyse the stability of solutions in terms of Floquet multipliers for Problem \ref{eq:nbody}, after rewriting it as the following first-order system of ODEs
\begin{equation}
    \begin{cases}
\label{e:first_ode}
    \dot{x} = M^{-1}\, y, \\
    \dot{y} = \nabla U(x),
\end{cases}
\end{equation}
where \( M^{-1} \in \mathbb R^{d\times n} \) is the inverse of the diagonal block-matrix 
\[
M^{-1}
 = 
 \begin{bmatrix}
     m_1\, I_n & 0 & 0
     \\
     0 & m_2\, I_n &0
     \\
     &  \ddots &
     \\
     0 & & m_d\, I_n
 \end{bmatrix}\]

In this case, \eqref{eq:matrix_eq} should be read with $A_{\bar x}$ the following block-matrix $\mathbb R^{(d\times n)^2}$:
\[
A_{\bar x} (t) = \begin{bmatrix}
    0 & M^{-1}
    \\
    \nabla^2 U(\bar x(t)) & 0
\end{bmatrix}
\]
where $0 =0_{\mathbb R^{d \times n}}$ and $\nabla^2 U$ is the Hessian matrix of the potential $U$ computed along the periodic solution $\bar x(t)$.

In the context of mechanical systems, the monodromy matrix \( \mathcal{M}_{\bar{x}} \) possesses the property that $\mathrm{det}(\mathcal{M}_{\bar{x}}) = 1$;  this implies that for every (complex) eigenvalue \( z \) with modulus \( |z| = r \), there exists another eigenvalue \( w \) such that \( |w| = 1/r \). As a result, an orbit \( \bar{x} \) is stable if and only if the spectrum of \( \mathcal{M}_{\bar x} \) is entirely contained within the unit circle \( \{|z| = 1\} \) in the complex plane.

 The following algorithm is proposed for reckoning the discrete Floquet stability of an orbit. 
\begin{algorithm}[H]
    \caption{Floquet Algorithm}
    \label{Floquet Algorithm}
    \begin{algorithmic}[1]
    \State \textbf{Input:} The periodic orbit, expressed in terms of its (truncated) Fourier coefficients.
    \State \textbf{Output:} The eigenvalues of the monodromy matrix. 
        \State \textbf{Step 1.} Build the Hessian matrix $\nabla^2U (\bar x (t))$
        \State \textbf{Step 2.} Define the matrix differential equation \eqref{eq:matrix_eq}
        \State \textbf{Step 3.} Compute the principal solution matrix $X(t)$ 
        \State \textbf{Step 4.} Compute the complex eigenvalues of $\mathcal M_{\bar x}=X(T)$
    \end{algorithmic}
\end{algorithm}
 
In Section \ref{sec:numeric}, we will apply the above algorithm to ten test cases, corresponding to orbits in the plane or the three-dimensional space with a variable number of bodies.

\subsection{Discrete Morse Index}

The second stability indicator we propose is the \textit{Morse index}, which measures the number of independent directions along which the action functional decreases, thus distinguishing minima, saddles, and unstable configurations within a variational framework~\cite{Morse1934, Smale1965}.

\noindent For numerical studies of periodic orbits and choreographies, the Morse index serves as a practical diagnostic tool, allowing one to detect bifurcations and qualitative changes in stability~\cite{Fukuda2019, Barutello2016}. Its discrete formulation, following the combinatorial approach of Forman~\cite{Forman1998}, adapts these ideas to finite-dimensional settings, where the action is defined on a discretised configuration space.


\noindent Following \cite{Ciccarelli2025}, we can give the following general definition.

\begin{definition}[Computational discrete Morse index]
\label{def 2: Discrete Morse index}
    Let \( f: H \to \mathbb{R} \) be a \( \mathcal{C}^2 \) function and \( p \) a non-degenerate critical point of \( f \). The \textit{index} of \( p \) is defined as the dimension of the maximal subspace of the tangent space at \( p \) where the Hessian is negative definite. Consequently, the \textit{discrete Morse index} \( \tilde{n}^{-} \) of \( p \) is the number of negative eigenvalues of the Hessian matrix \( H_f(p) \).
\end{definition}

\noindent In summary, the discrete Morse index provides a numerical way to distinguish between different types of critical points: minima satisfy \( \tilde{n}^{-}(p) = 0 \), while saddle points correspond to \( \tilde{n}^{-}(p) > 0 \).

\subsubsection{Application to the $n$-body problem} 
In the context of the $n$-body problem, the computation of the discrete Morse index requires first the discretisation of the action functional. This preliminary step transforms the continuous variational problem into a finite-dimensional optimisation problem, whose critical points approximate periodic orbits of the system. Once the action functional has been discretised, the Morse index can then be evaluated either on the \textit{fundamental domain} $\I$ or, by exploiting the symmetries of the orbit, on the \textit{entire domain}. Clearly, the index computed on the fundamental domain does not exceed that obtained over the full orbit.

\paragraph{Discretisation of the action functional}
\noindent Following \cite{symorb}, the action functional \eqref{action functional} can be discretised in two main ways.
\paragraph{Point discretisation}
\label{Discretisation through points}
In this approach, the integral is approximated by quadrature on a uniform grid 
\( 0 = t_0 < t_1 < \dots < t_{M+1} = {T} \). 
Derivatives are replaced with finite differences, and the action is written in terms of the discrete variables 
\( y_i^k \approx x_i(t_k) \in \mathbb{R}^d \), for $i=1,\dots,n$ and $k=0,\dots,M+1$. 
This yields the discrete functional
\begin{equation}
\label{Approximated action functional}
    f_1(y_i^k) = \mathcal{A}^{(1)}_h,
\end{equation}
defined on the block matrix \( (y_i^k)_{i,k} \in \mathbb{R}^{d\times n \times (M+2)} \). 
For details see \cite{CI24, symorb}.
\paragraph{Fourier coefficients discretisation} Alternatively, the trajectory is approximated by a truncated Fourier series added to a linear interpolation between the boundary configurations:
\begin{equation}
\label{eq:F_tronc}
    x(t) = x_0 + \tfrac{t}{\pi}(x_1-x_0) + \sum_{k=1}^F A_k \sin(kt), \quad t \in [0,\pi],
\end{equation}
where the sine terms ensure periodicity. Substituting into \eqref{action functional} gives
\begin{equation}
\label{eq:f2}
    f_2(x_0,x_1,A_k) = \mathcal{A}^{(2)}_h.
\end{equation}
The variables are the endpoints and Fourier coefficients, with dimension \((F+2)\times d \times n\). 
See \cite{symorb, Sim, Sim2} for further discussion. \\

\noindent In both schemes, the problem reduces to optimising an objective function \( f = f_i \) (\(i=1,2\)) with respect to the chosen discrete variables, whose critical points correspond to periodic solutions of the system.
Once a discrete periodic orbit has been obtained, one can further characterise the nature of the corresponding critical point by computing its discrete Morse index.
This can be done in two equivalent settings, depending on whether one exploits the symmetries of the solution or not: either on the fundamental domain, where only a representative segment of the orbit is considered, or on the full orbit, where the computation is performed over the entire period.

\paragraph{\textit{On the fundamental domain}}  
As described in \cite{symorbJL, Ciccarelli2025}, the Hessian of the discretised action functional can be computed within the fundamental domain \( \mathbb{I} \) either by using a Fourier-based discretisation, where the trajectory is approximated by truncated series, or by direct pointwise discretisation, where the time interval is subdivided into a finite number of nodes.

\paragraph{\textit{On the full orbit}}  
Alternatively, the discrete Morse index can be computed on the full periodic orbit. In this case, the action functional is discretised over the entire period, ensuring that periodic boundary conditions are satisfied. Since a periodic orbit has identical initial and final points, only one of them is retained in the computation.

\subsubsection{Algorithm}

The following procedure summarises the computation of the discrete Morse index on the entire orbit.

\begin{algorithm}[H]
    \caption{Computation of the Discrete Morse Index}
    \label{alg:morse_index}
    \begin{algorithmic}[1]
    \State \textbf{Input:} Periodic solution given either 
        \begin{itemize}
            \item[(a)] by its Fourier coefficients, or 
            \item[(b)] by sampled positions $x_i(t_k)$ for each body $i=1,\dots,n$ at times $t_k$, $k=1,\dots,w$.
        \end{itemize}
    \State \textbf{Output:} Discrete Morse index $\tilde{n}^{-}(p)$.
        \State \textbf{Step 1.} Build the Hessian matrix of the discretised action functional:
        \begin{itemize}
            \item Fourier case: compute the Hessian in coefficient space.  
            \item Pointwise case: assemble the Hessian from the discretised action.
        \end{itemize}
        \State \textbf{Step 2.} Compute all eigenvalues $\lambda_j$ of the Hessian.
        \State \textbf{Step 3.} The discrete Morse index is the number of negative eigenvalues:
        \[
        \tilde{n}^{-}(p) = \#\{\lambda_j < 0\}.
        \]
    \end{algorithmic}
\end{algorithm}

\noindent It is worth noting that, in contrast to previous studies that analysed specific families of solutions (such as the figure-eight choreography \cite{chen_mont00}), the approach adopted here is formulated to be applicable to any periodic orbit obtained numerically. In this way, the abstract concepts introduced earlier are translated into general computational tools, establishing a direct connection between the theoretical framework and the analysis of stability in the $n$-body problem.

\section{Numerical results} \label{sec:numeric}
\begin{figure}[htbp!]
\centering
\makebox[\textwidth][c]{
\fbox{%
\begin{minipage}{1\linewidth}
\centering

\begin{tabular}{|m{4cm}|l|l|}
\hline
\textbf{Group generator(s)}  & \multicolumn{2}{|l|}{
\begin{tabular}{@{}c@{}}
$\ker \tau = \langle \Id \rangle, \quad \bar G = \langle r \rangle $ \\ 
$\rho(r) = - \Id, \quad\sigma(r) = ()$
\end{tabular}} \\
\hline
\multicolumn{2}{|l|}{\textbf{Action value}} & 9.8022\\
\hline
\multicolumn{2}{|l|}{\textbf{Gradient norm}} & $7.72\times 10^{-11}$\\
\hline
\multirow{2}{10em}{\textbf{Morse (fund. domain)}} & \textbf{Index} & 0\\
\cline{2-3}
 & \textbf{Max negative eigenvalue} & ---\\
\hline
\multirow{2}{10em}{\textbf{Morse (period)}} & \textbf{Index} & 0\\
\cline{2-3}
 & \textbf{Max negative eigenvalue} & ---\\
\hline
\textbf{Floquet} 
 & \textbf{Norm of the max eigenvalue} & 85.0213\\
\hline
\end{tabular}

\vspace{1em}

\includegraphics[width=0.8\linewidth]{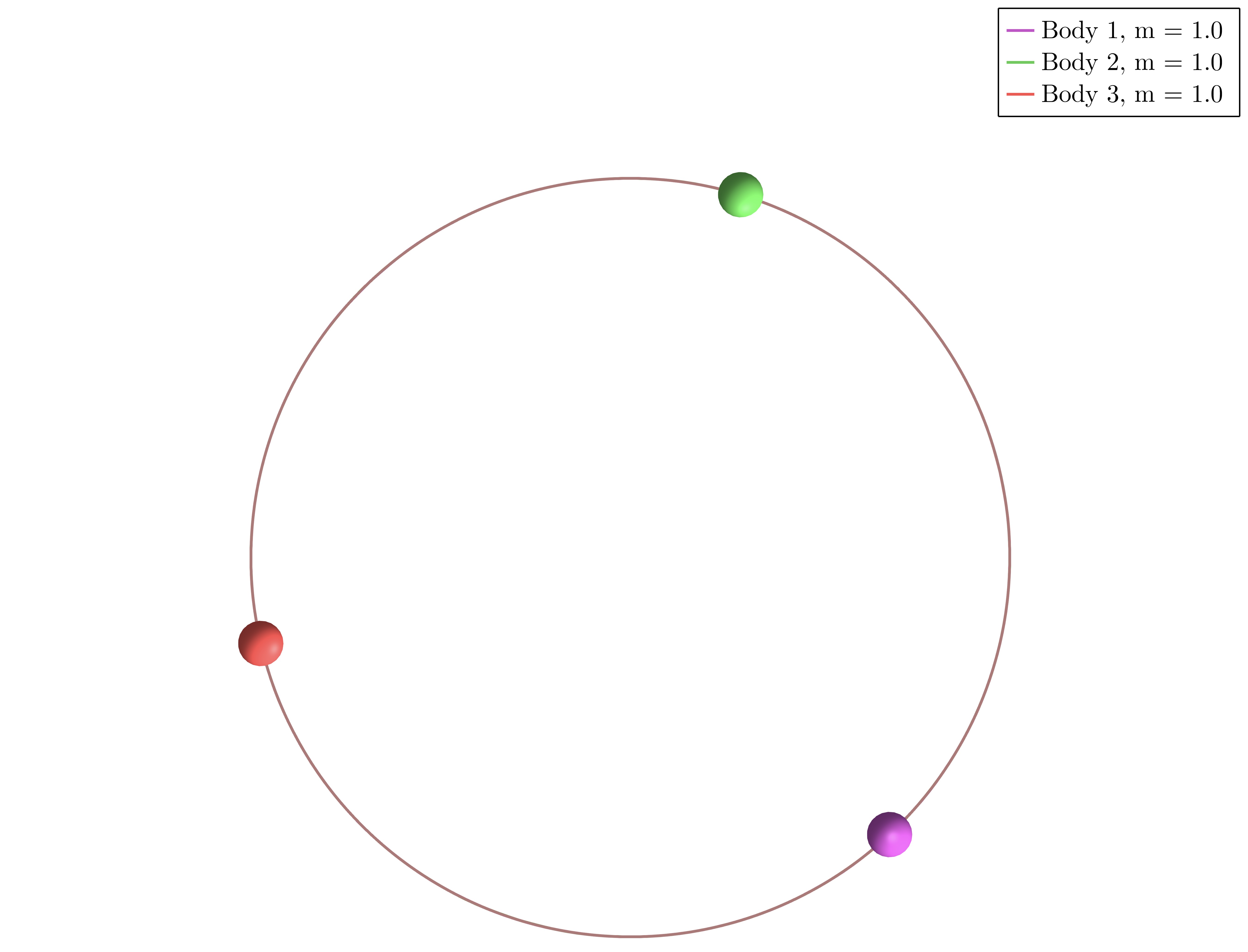} 
\end{minipage}
}}
\caption{}
\label{fig:fig1}
\end{figure}

\begin{figure}[htbp!]
\centering
\makebox[\textwidth][c]{
\fbox{%
\begin{minipage}{1\linewidth}
\centering

\begin{tabular}{|m{4cm}|l|l|}
\hline
\textbf{Group generator(s)}  & \multicolumn{2}{|l|}{
\begin{tabular}{@{}c@{}}
$\ker \tau = \langle \Id \rangle, \quad \bar G = \langle r \rangle $ \\ 
$\rho(r) = - \Id, \quad\sigma(r) = ()$
\end{tabular}
} \\
\hline
\multicolumn{2}{|l|}{\textbf{Action value}} & 10.9365\\
\hline
\multicolumn{2}{|l|}{\textbf{Gradient norm}} & $6.77\times 10^{-11}$\\
\hline
\multirow{2}{10em}{\textbf{Morse (fund. domain)}} & \textbf{Index} & 1\\
\cline{2-3}
 & \textbf{Max negative eigenvalue} & -5.5278\\
\hline
\multirow{2}{10em}{\textbf{Morse (period)}} & \textbf{Index} & 3\\
\cline{2-3}
 & \textbf{Max negative eigenvalue} &  -13.9510\\
\hline
\textbf{Floquet} 
 & \textbf{Norm of the max eigenvalue} & $5.8982\times 10^4$\\
\hline
\end{tabular}

\vspace{1em}

\includegraphics[width=0.8\linewidth]{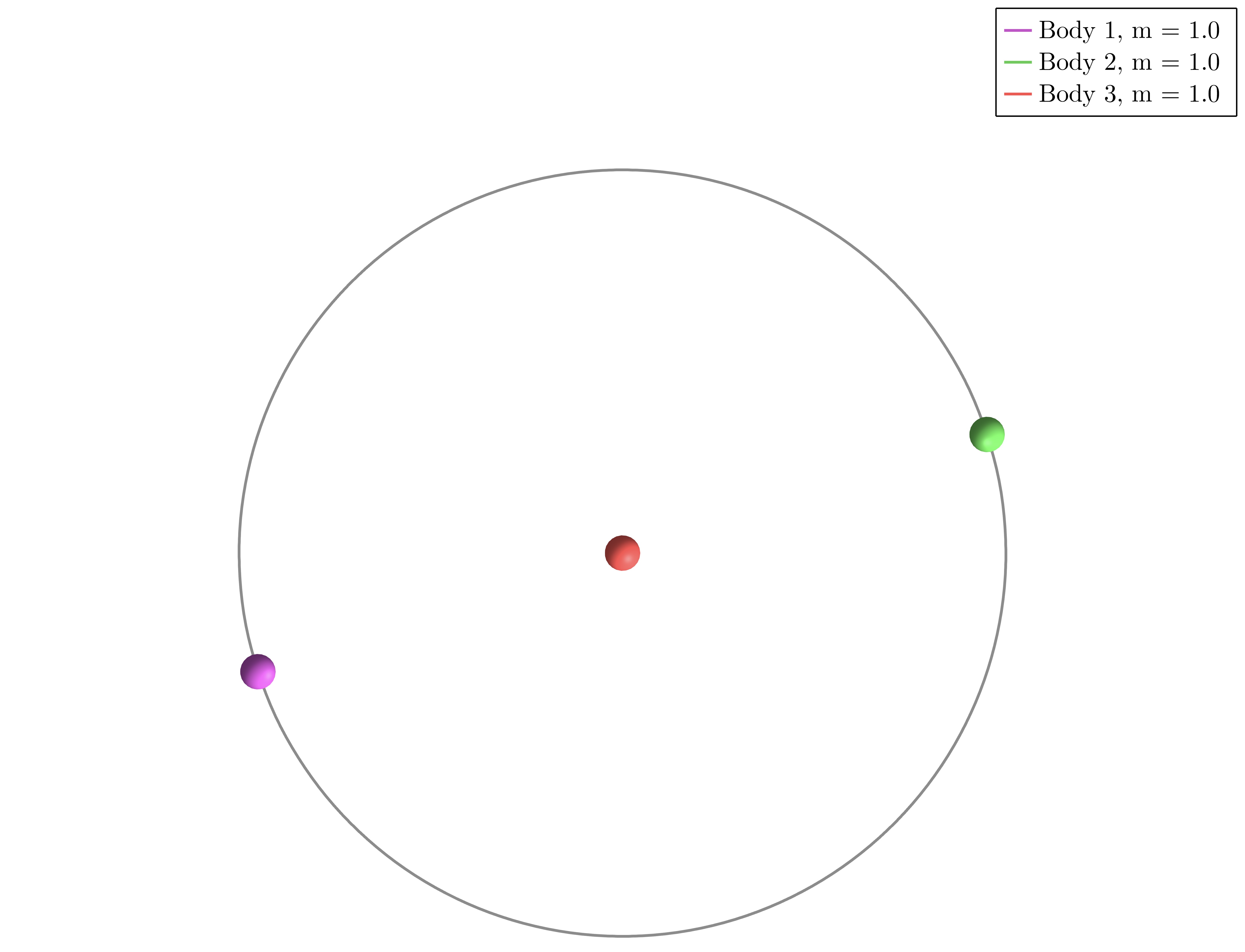} 

\end{minipage}
}}
\caption{}
\label{fig:fig2}
\end{figure}

\begin{figure}[htbp!]
\centering
\makebox[\textwidth][c]{
\fbox{%
\begin{minipage}{1\linewidth}
\centering

\begin{tabular}{|m{4cm}|l|l|}
\hline
\textbf{Group generator(s)}  & \multicolumn{2}{|l|}{
\begin{tabular}{@{}l@{}}
$\ker \tau = \langle \Id \rangle, \quad \bar G = \langle r, s \rangle $ \\ 
$\rho(r) = \Id, \quad\sigma(r) = (1,2,3)$ \\ 
$\rho(s) = -\Id, \quad\sigma(s) = (2,3)$ 
\end{tabular}
} \\
\hline
\multicolumn{2}{|l|}{\textbf{Action value}} & 5.8584\\
\hline
\multicolumn{2}{|l|}{\textbf{Gradient norm}} & $4.85\times 10^{-12}$\\
\hline
\multirow{2}{10em}{\textbf{Morse (fund. domain)}} & \textbf{Index} & 0\\
\cline{2-3}
 & \textbf{Max negative eigenvalue} & ---\\
\hline
\multirow{2}{10em}{\textbf{Morse (period)}} & \textbf{Index} & 2\\
\cline{2-3}
 & \textbf{Max negative eigenvalue} & -0.2257\\
\hline
\textbf{Floquet}
 & \textbf{Norm of the max eigenvalue} & 1.0187\\
\hline
\end{tabular}

\vspace{1em}

\includegraphics[width=0.8\linewidth]{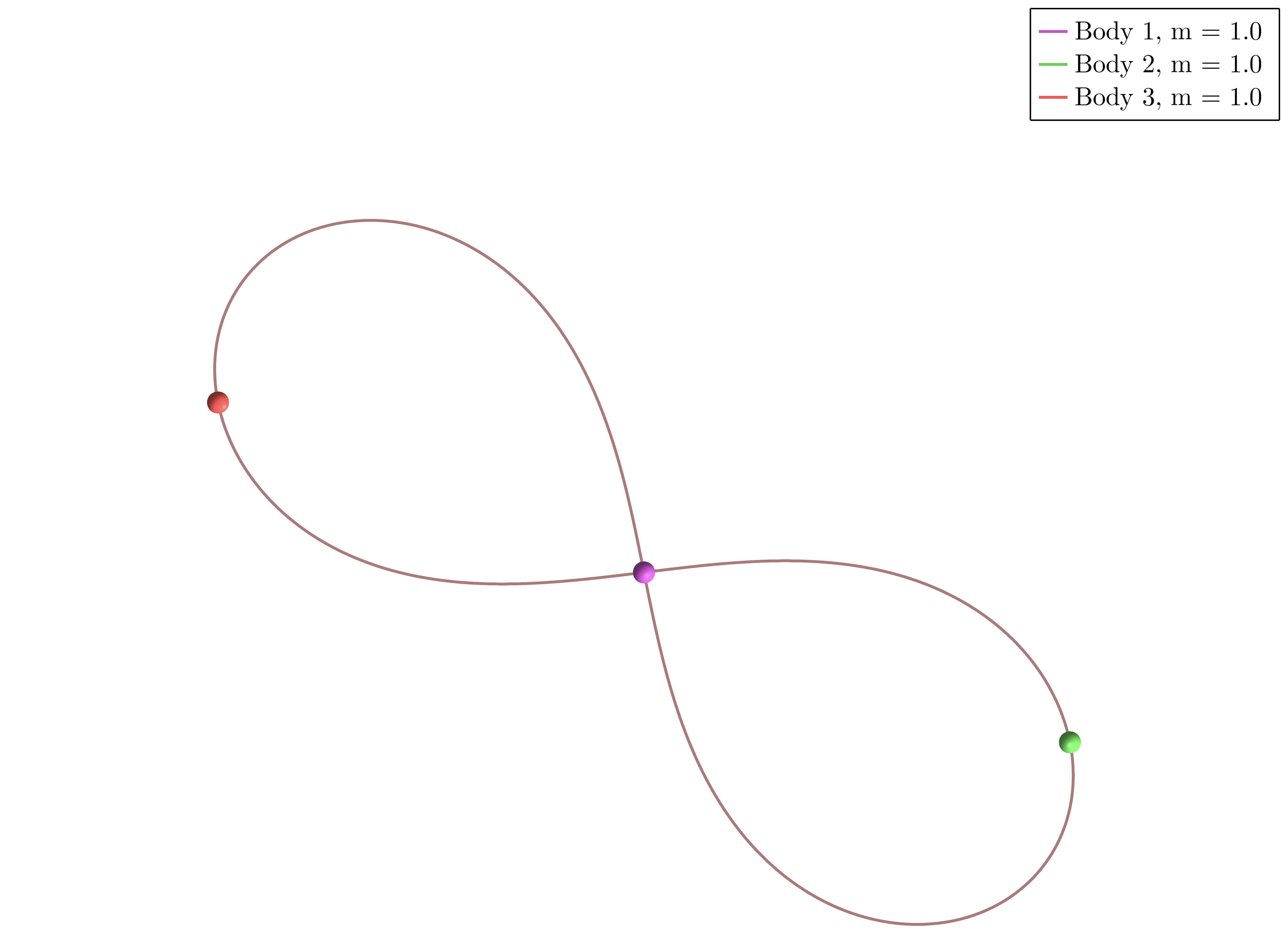} 

\end{minipage}
}}
\caption{}
\label{fig:fig3}
\end{figure}

\begin{figure}[htbp!]
\centering
\makebox[\textwidth][c]{
\fbox{%
\begin{minipage}{1\linewidth}
\centering

\begin{tabular}{|m{4cm}|l|l|}
\hline
\textbf{Group generator(s)}  & \multicolumn{2}{|l|}{
\begin{tabular}{@{}l@{}}
$\ker \tau = \langle \Id \rangle, \quad \bar G = \langle r \rangle $ \\ 
$\rho(r) = -\Id, \quad\sigma(r) = ()$ \\ 
\end{tabular}
} \\
\hline
\multicolumn{2}{|l|}{\textbf{Action value}} & 10.4421\\
\hline
\multicolumn{2}{|l|}{\textbf{Gradient norm}} & $3.22\times 10^{-8}$\\
\hline
\multirow{2}{10em}{\textbf{Morse (fund. domain)}} & \textbf{Index} & 0\\
\cline{2-3}
 & \textbf{Max negative eigenvalue} & ---\\
\hline
\multirow{2}{10em}{\textbf{Morse (period)}} & \textbf{Index} & 2\\
\cline{2-3}
 & \textbf{Max negative eigenvalue} & -0.0617\\
\hline
\textbf{Floquet}
 & \textbf{Norm of the max eigenvalue} & 1.0462\\
\hline
\end{tabular}

\vspace{1em}

\includegraphics[width=0.8\linewidth]{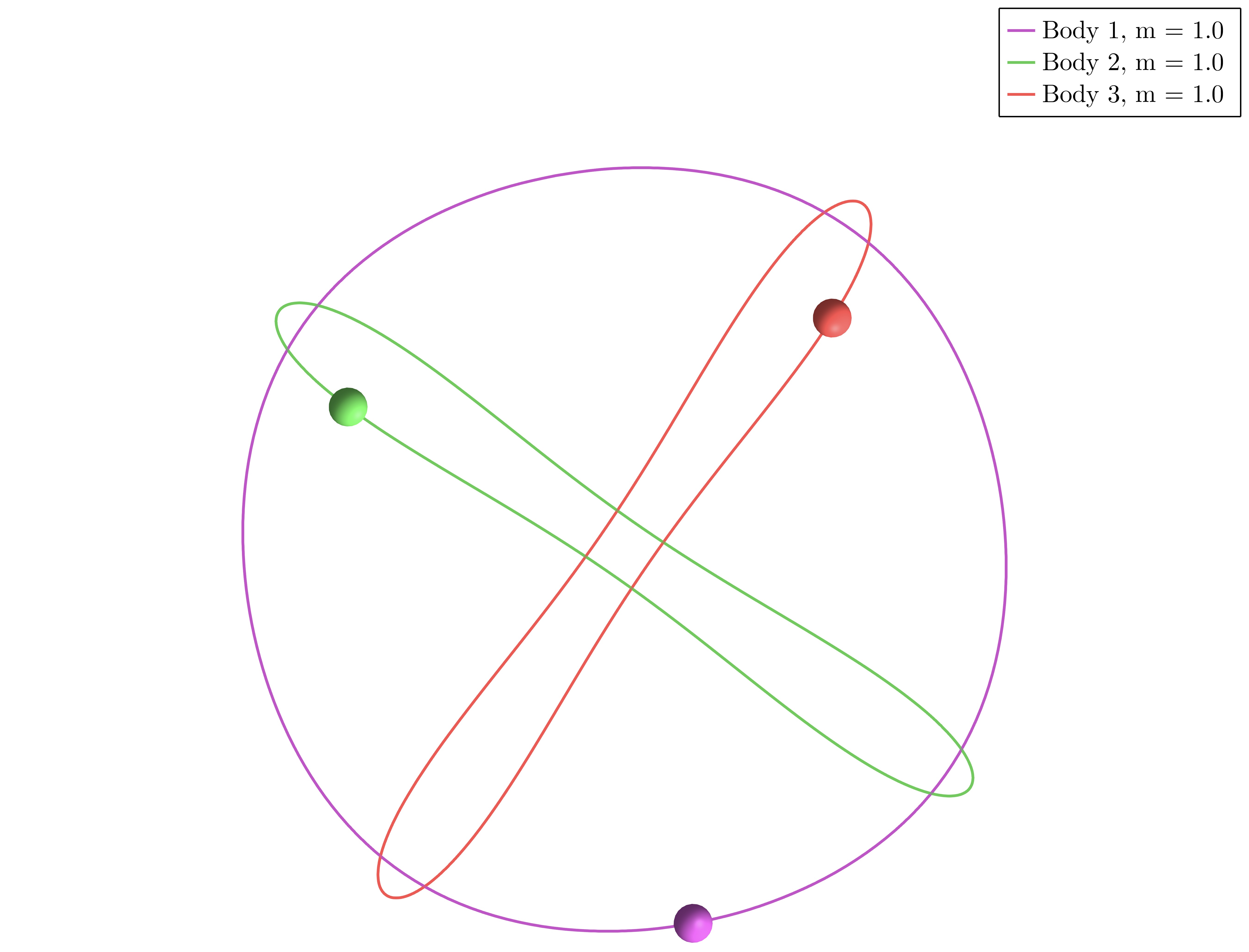} 

\end{minipage}
}}

\caption{}
\label{fig:fig4}
\end{figure}

\begin{figure}[htbp!]
\centering
\makebox[\textwidth][c]{
\fbox{%
\begin{minipage}{1\linewidth}
\centering

\begin{tabular}{|m{4cm}|l|l|}
\hline
\textbf{Group generator(s)}  & \multicolumn{2}{|l|}{
\begin{tabular}{@{}l@{}}
$\ker \tau = \langle \Id \rangle, \quad \bar G = \langle r \rangle $ \\ 
$\rho(r) = \Id, \quad\sigma(r) = (1,2,3,4,5,6,7,8,9,10,11,12)$ \\ 
 
\end{tabular}
} \\
\hline
\multicolumn{2}{|l|}{\textbf{Action value}} & 72.6212\\
\hline
\multicolumn{2}{|l|}{\textbf{Gradient norm}} & $2.01\times10^{-8}$\\
\hline
\multirow{2}{10em}{\textbf{Morse (fund. domain)}} &\textbf{Index} & 0\\
\cline{2-3}
 & \textbf{Max negative eigenvalue} & ---\\
\hline
\multirow{2}{10em}{\textbf{Morse (period)}} & \textbf{Index} & 4\\
\cline{2-3}
 & \textbf{Max negative eigenvalue} & -0.1103\\
\hline
\textbf{Floquet} 
 & \textbf{Norm of the max eigenvalue} & $4.1930\times 10^{12}$\\
\hline
\end{tabular}

\vspace{1em}

\includegraphics[width=0.7\linewidth]{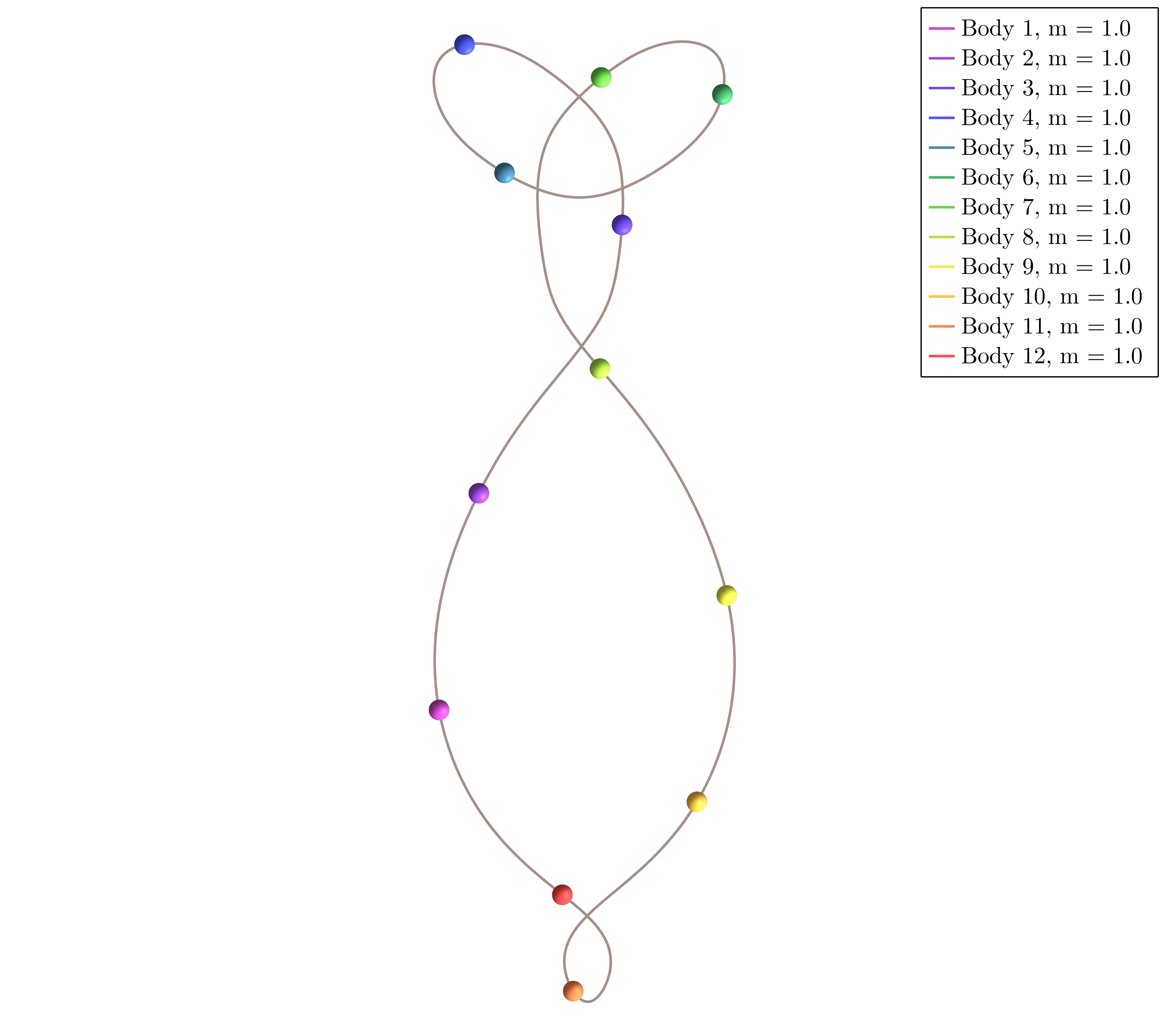} 

\end{minipage}
}}
\caption{}
\label{fig:fig5}
\end{figure}

\begin{figure}[htbp!]
\centering
\makebox[\textwidth][c]{
\fbox{%
\begin{minipage}{1\linewidth}
\centering

\begin{tabular}{|m{4cm}|l|l|}
\hline
\textbf{Group generator(s)}  & \multicolumn{2}{|l|}{
\begin{tabular}{@{}l@{}}
$\ker \tau = \langle \kappa \rangle, \quad \bar G = \langle r \rangle $ \\ 
$\rho(\kappa) = R(2\pi/5)$, $\sigma(\kappa) = (1,2,3,4,5)(6,7,8,9,10)$ \\ 
$\rho(r) = -\Id, \quad\sigma(r) = ()$ 
\end{tabular}
} \\
\hline
\multicolumn{2}{|l|}{\textbf{Action value}} & 114.3331\\
\hline
\multicolumn{2}{|l|}{\textbf{Gradient norm}} & $4.11\times 10^{-7}$\\
\hline
\multirow{2}{10em}{\textbf{Morse (fund. domain)}} & \textbf{Index} & 0\\
\cline{2-3}
 & \textbf{Max negative eigenvalue} & ---\\
\hline
\multirow{2}{10em}{\textbf{Morse (period)}} & \textbf{Index }& 2\\
\cline{2-3}
 & \textbf{Max negative eigenvalue} & -0.0024\\
\hline
\textbf{Floquet} 
 & \textbf{Norm of the max eigenvalue} &  1173.9612\\
\hline
\end{tabular}

\vspace{1em}

\includegraphics[width=0.7\linewidth]{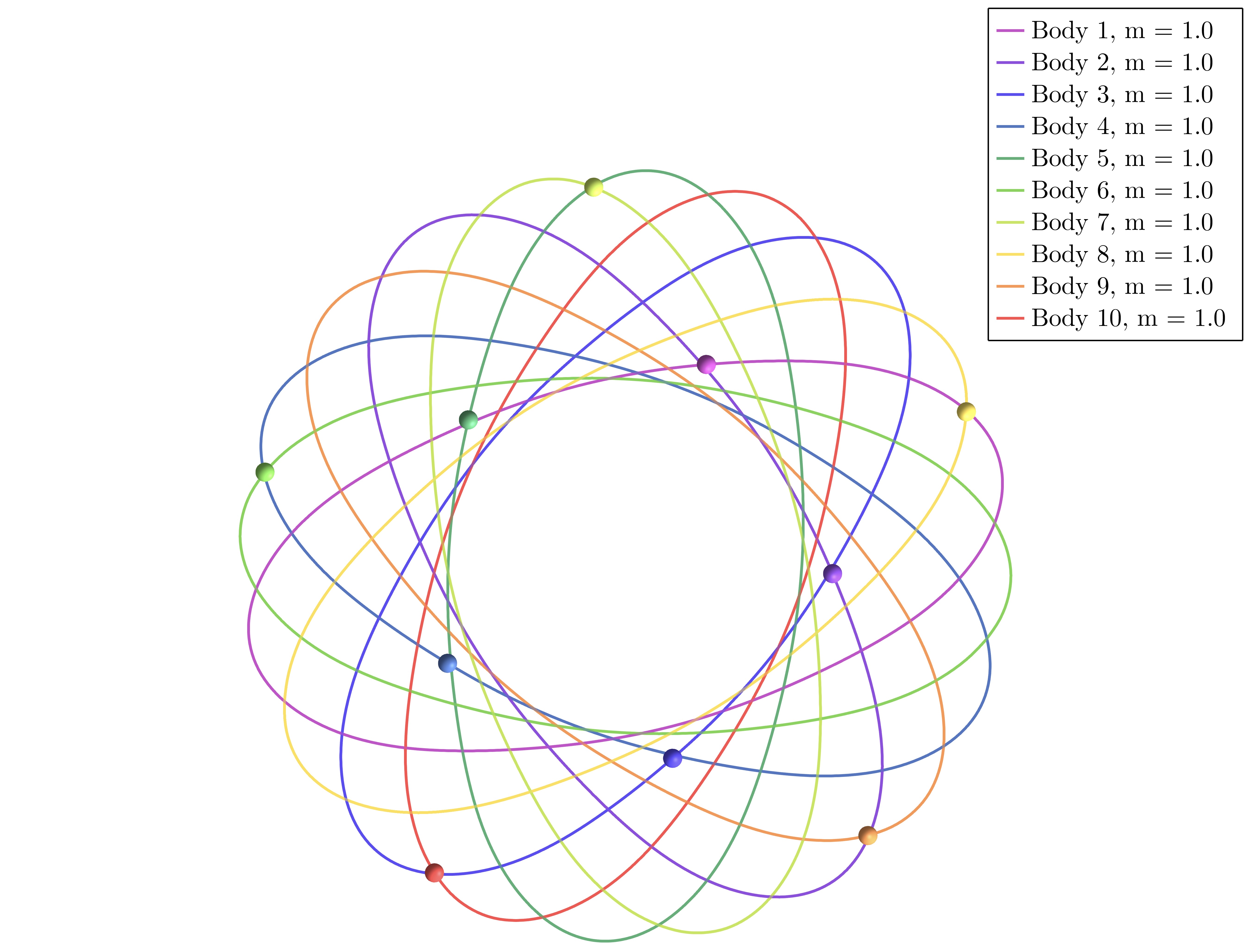} 

\end{minipage}
}}

\caption{}
\label{fig:fig6}
\end{figure}

\begin{figure}[htbp!]
\centering
\makebox[\textwidth][c]{
\fbox{%
\begin{minipage}{1\linewidth}
\centering

\begin{tabular}{|m{4cm}|l|l|}
\hline
\textbf{Group generator(s)}  & \multicolumn{2}{|l|}{
\begin{tabular}{@{}l@{}}
$\ker \tau = \langle \Id \rangle, \quad \bar G = \langle r \rangle $ \\ 
$\rho(r) =\small{\begin{pmatrix}
    R(2\pi/3) & 0\\ 0 & -1
\end{pmatrix}}, \quad \sigma(r) = ()$ \\ 
\end{tabular}
} \\
\hline
\multicolumn{2}{|l|}{\textbf{Action value}} & 8.3409\\
\hline
\multicolumn{2}{|l|}{\textbf{Gradient norm}} & $3.46\times 10^{-8}$\\
\hline
\multirow{2}{10em}{\textbf{Morse (fund. domain)}} & \textbf{Index }& 1\\
\cline{2-3}
 & \textbf{Max negative eigenvalue} & -3.0423\\
\hline
\multirow{2}{10em}{\textbf{Morse (period)}} &\textbf{Index} & 19\\
\cline{2-3}
 & \textbf{Max negative eigenvalue} & -48.5411\\
\hline
\textbf{Floquet} 
 & \textbf{Norm of the max eigenvalue} & $6.3110\times 10^8$\\
\hline
\end{tabular}

\vspace{1em}

\includegraphics[width=0.8\linewidth]{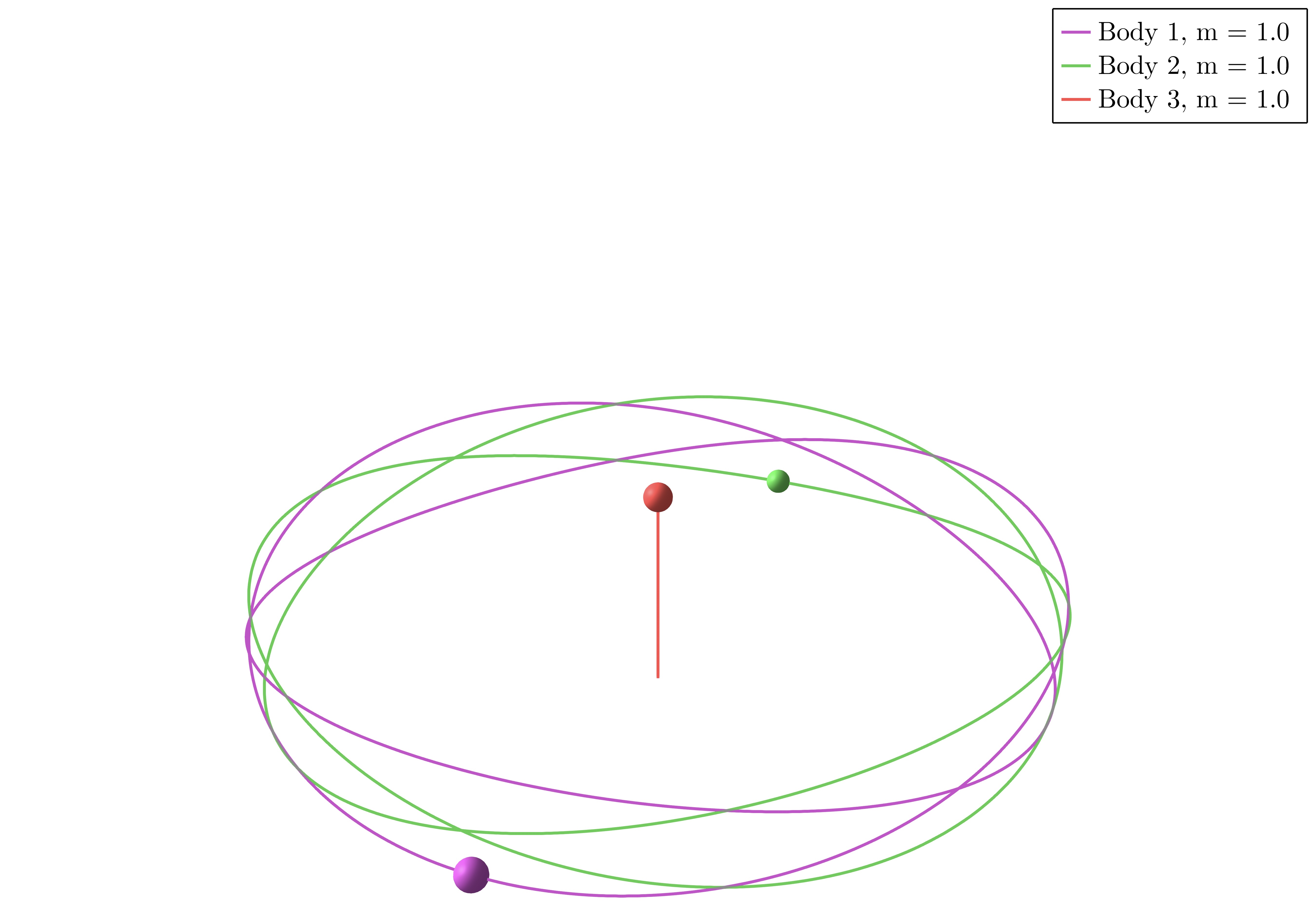} 

\end{minipage}
}}
\caption{}
\label{fig:fig7}
\end{figure}

\begin{figure}[htbp!]
\centering
\makebox[\textwidth][c]{
\fbox{%
\begin{minipage}{1\linewidth}
\centering

\begin{tabular}{|m{4cm}|l|l|}
\hline
\textbf{Group generator(s)} & \multicolumn{2}{|l|}{
\begin{tabular}{@{}l@{}}
$\ker \tau = \langle \Id \rangle, \quad \bar G = \langle r, s \rangle $ \\ 
$\rho(r) = - \Id, \quad \sigma(r) = ()$ \\ 
$\rho(s) = \small{\begin{pmatrix}
-Id_2 & 0\\ 0 &1
\end{pmatrix}}, \quad \sigma(s) = ()$ \\ 
\end{tabular}
} \\
\hline
\multicolumn{2}{|l|}{\textbf{Action value}} & 11.8999\\
\hline
\multicolumn{2}{|l|}{\textbf{Gradient norm}} & $2.52\times 10^{-8}$\\
\hline
\multirow{2}{10em}{\textbf{Morse (fund. domain)}} & \textbf{Index }& 0\\
\cline{2-3}
 & \textbf{Max negative eigenvalue} & ---\\
\hline
\multirow{2}{10em}{\textbf{Morse (period)}} & \textbf{Index} & 9\\
\cline{2-3}
 &\textbf{Max negative eigenvalue} & -11.6320\\
\hline
\textbf{Floquet} 
 &\textbf{Norm of the max eigenvalue}& 25.1735\\
\hline
\end{tabular}

\vspace{1em}

\includegraphics[width=0.7\linewidth]{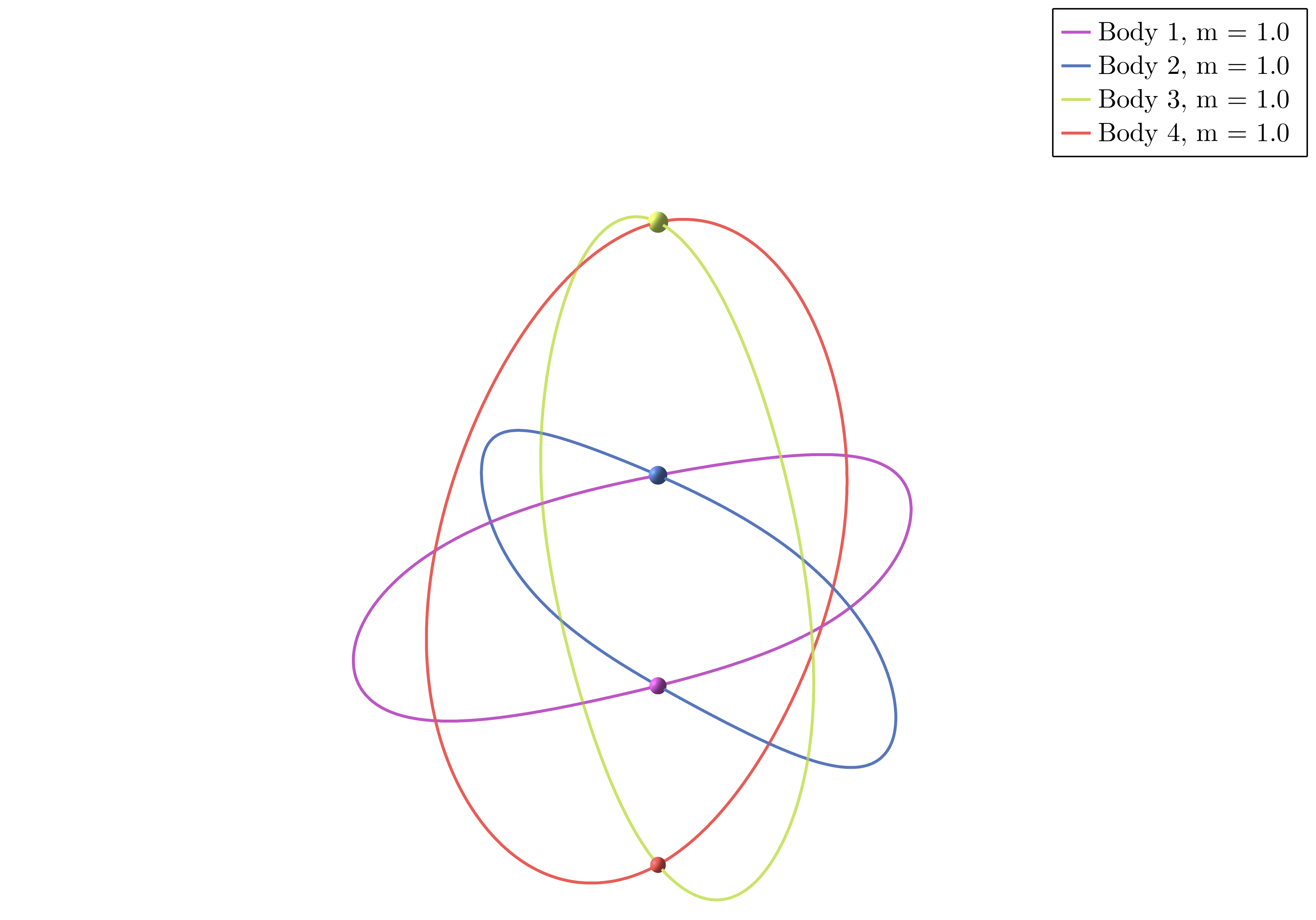}

\end{minipage}
}}

\caption{}
\label{fig:fig8}
\end{figure}

\begin{figure}[htbp!]
\centering
\makebox[\textwidth][c]{
\fbox{%
\begin{minipage}{1\linewidth}
\centering

\begin{tabular}{|m{4cm}|l|l|}
\hline
\textbf{Group generator(s)}  & \multicolumn{2}{|l|}{
\begin{tabular}{@{}l@{}}
$\ker \tau = \langle \Id \rangle, \quad \bar G = \langle r, s \rangle $ \\ 
$\rho(r) = - \Id \quad \sigma(r) = (1,2)(3,4)$\\
$\rho(s) = \small{\begin{pmatrix}
    -\Id_2&0\\0&1
\end{pmatrix}} \quad \sigma(s) = ()$ \\ 
\end{tabular}
} \\
\hline
\multicolumn{2}{|l|}{\textbf{Action value}} & 12.1428\\
\hline
\multicolumn{2}{|l|}{\textbf{Gradient norm}} & $5.27\times 10^{-10}$\\
\hline
\multirow{2}{10em}{\textbf{Morse (fund. domain)}} & \textbf{Index }& 0\\
\cline{2-3}
 & \textbf{Max negative eigenvalue} & ---\\
\hline
\multirow{2}{10em}{\textbf{Morse (period)}} &\textbf{Index }& 11\\
\cline{2-3}
 & \textbf{Max negative eigenvalue} & -202.1842\\
\hline
\textbf{Floquet} 
 & \textbf{Norm of the max eigenvalue } & 2191.5747\\
\hline
\end{tabular}

\vspace{1em}

\includegraphics[width=0.65\linewidth]{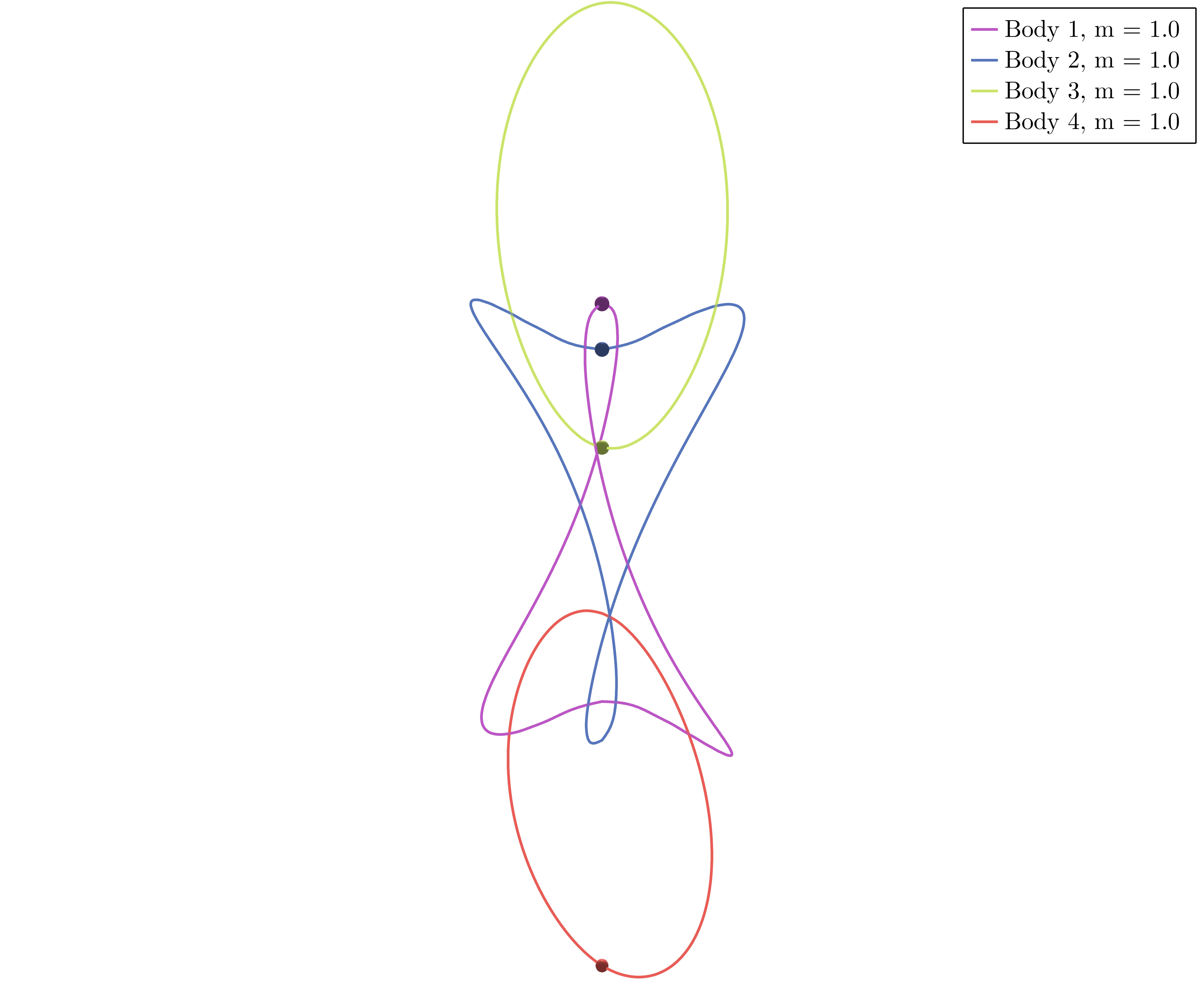} 

\end{minipage}
}}
\caption{}
\label{fig:fig9}
\end{figure}

\begin{figure}[htbp!]
\centering
\makebox[\textwidth][c]{
\fbox{%
\begin{minipage}{1\linewidth}
\centering

\begin{tabular}{|m{4cm}|l|l|}
\hline
\textbf{Group generator(s)}  & \multicolumn{2}{|l|}{
\begin{tabular}{@{}l@{}}
$\ker \tau = \langle \Id \rangle, \quad \bar G = \langle r \rangle $ \\ 
$\rho(r) = R(\pi/6)$\\$\sigma(r) = (1,2,3,4,5,6,7,8,9,10,11,12)$ \\ 
\end{tabular}
} \\
\hline
\multicolumn{2}{|l|}{\textbf{Action value}} & 81.6410\\
\hline
\multicolumn{2}{|l|}{\textbf{Gradient norm}} & $7.65\times 10^{-9}$\\
\hline
\multirow{2}{10em}{\textbf{Morse (fund. domain)}} &\textbf{Index}& 0\\
\cline{2-3}
 & \textbf{Max negative eigenvalue} & ---\\
\hline
\multirow{2}{10em}{\textbf{Morse (period)}} & \textbf{Index} & 4\\
\cline{2-3}
 & \textbf{Max negative eigenvalue} & -0.1778\\
\hline
\textbf{Floquet} 
 & \textbf{Norm of the max eigenvalue }&  $8.5410\times 10^{12}$\\
\hline
\end{tabular}

\vspace{1em}

\includegraphics[width=0.7\linewidth]{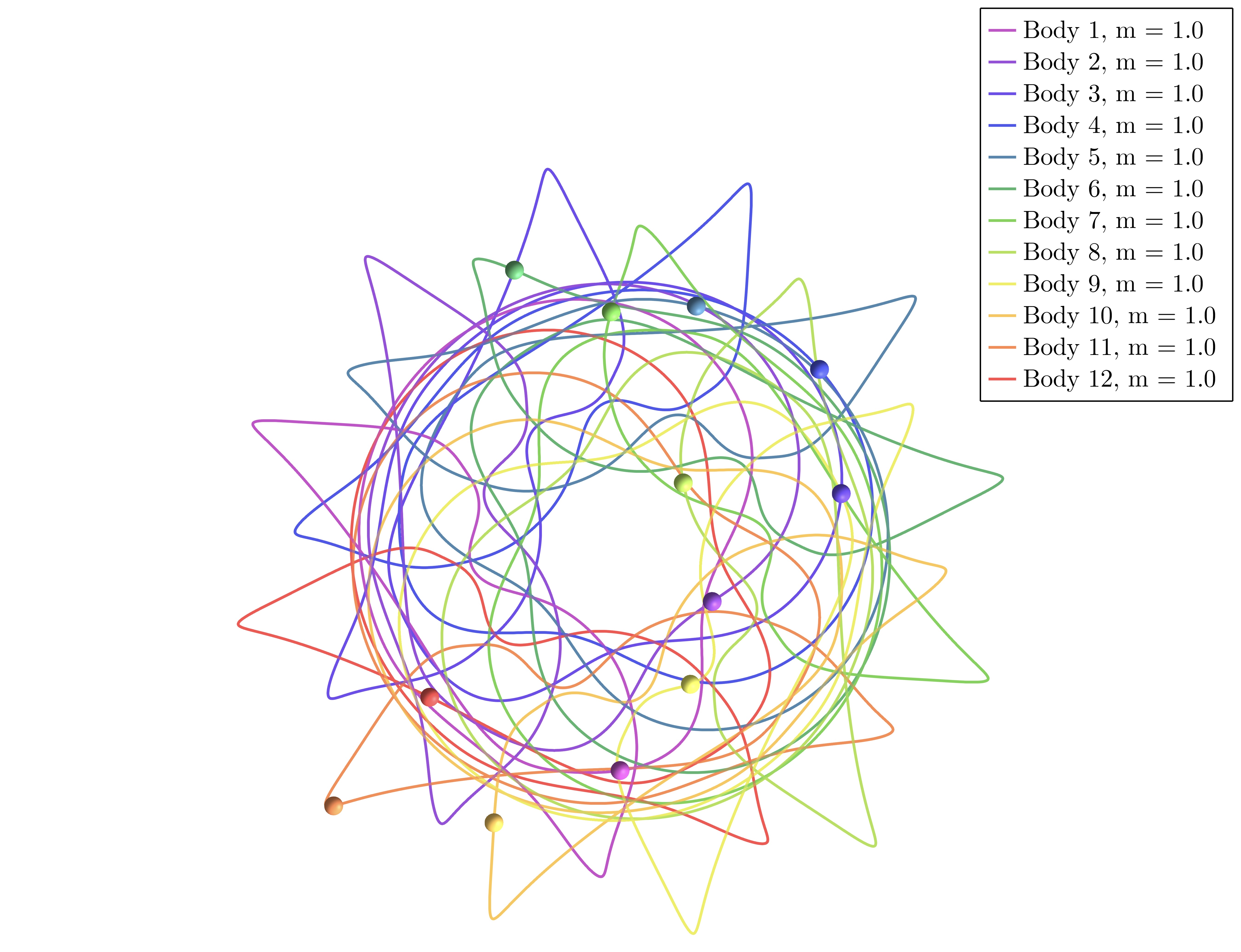} 

\end{minipage}
}}
\caption{}
\label{fig:fig10}
\end{figure}

This Section presents the application of the algorithms described in Section \ref{sec:stab} to a set of ten distinct periodic obits computed using \textit{SymOrb.jl}. Figures \ref{fig:fig1}-\ref{fig:fig10} display a wide variety of examples, including both planar and spatial configurations, and involving different numbers of bodies. For each orbit, we report:
\begin{itemize}
    \item the generators of the finite symmetry group $G$ associated with the orbits; recall that $\bar G = G/\ker\tau$ (see Section \ref{sec:intro}), and that the notation $R(\theta)$ denotes the rotation matrices of the form 
    \[
    R(\theta)=\begin{pmatrix}
        \cos\theta & -\sin\theta\\\sin\theta &\cos\theta
    \end{pmatrix};
    \]
    \item the values of the action $\mathcal{A}$ and its gradient, to assess the accuracy of the numerical method in locating a critical point of $\mathcal A$; 
    \item the numerical Morse index, computed both on the fundamental domain $\mathbb{I}$ and over the entire period; whenever this index is non-zero, we additionally report the maximal negative eigenvalue of $Hess(\mathcal A)$ ;
    \item the largest modulus of the eigenvalues of the monodromy matrix, i.e., the maximal Floquet multiplier. 
\end{itemize} 
In the variational framework of the Morse index, the maximal negative eigenvalue provides a quantitative indication on how clearly the negative spectrum is separated from zero, and thus on the reliability of the computed Morse index. Conversely, in the Floquet analysis, the focus lies on the maximal (in modulus) eigenvalue of the monodromy matrix: when it lies far from $\pm 1$, we can confidently conclude that the orbit, though numerical approximated, is linearly unstable. \\
The selected set of orbits is designed to represent the diversity of periodic solutions we can find in the $n-$body problem, with particular attention to non-collisional examples.\\
A notable subset consists in the so-called \textit{choreographies}, where all bodies move along a common trajectory, as in Figures \ref{fig:fig1}, \ref{fig:fig3}, and \ref{fig:fig5}. Moreover, Figures \ref{fig:fig1}-\ref{fig:fig4} correspond to classical solutions already studied in the literature (see, e.g., \cite{Barut16, Roberts2002, Simo, Roberts, vanderbei2004new}), included here to facilitate comparison between our computations and known results. In particular, Figure \ref{fig:fig2} shows a {\em collinear relative equilibrium}, a well-known solution whose first study traces back to early works such as \cite{Meyer1933}.\\
Among our findings, it is worth noticing that the triangular Lagrange solution (Figure \ref{fig:fig1}) minimises the action functional, despite being linearly unstable. Conversely, the celebrated figure-eight choreography (Figure \ref{fig:fig3}) is not a minimiser of the action functional over its full period, yet it exhibits substantial linear stability -- in agreement with previous studies (\cite{Simo, Roberts, galan, KapelaSimo}). \\
Overall, our results illustrate the rich interplay between the variational and dynamical features of periodic orbits in the $n$-body problem. The comparison between Morse index (reflecting the local behaviour of $\mathcal{A}$ as a variational functional) and Floquet multipliers (characterising linear stability) reveals a wide spectrum of behaviours, which may serve as a useful tool for classifying and distinguishing periodic solutions according to their shared dynamical properties.

\newpage

\paragraph{\textbf{Funding}} 
M.I. conducted this research during and with the support of the Italian national inter-university PhD programme in Space Science and Technology. Moreover M.I. is partially supported by INdAM group G.N.A.M.P.A.\\
I.D.B. is supported by the Italian Research Center on High Performance Computing Big Data and Quantum Computing (ICSC), project funded by European Union - NextGenerationUE and National Recovery and Resilience Plan (NRRP) - Mission 4 Component 2. Spoke 3, Astrophysics and Cosmos Observations. D. B. was supported by the project NODES which has received funding from the MUR – M4C2 1.5 of PNRR funded by the European Union – NextGenerationEU (Grant agreement no. ECS00000036). D. B. was partially supported by INdAM group G.N.A.M.P.A.
 D.P. is partially supported by INdAM group G.N.A.M.P.A.

\paragraph{\textbf{Declaration of generative AI in scientific writing}} During the revision of this work the authors used LLMs to the sole scope of review the language. After using this tools, the authors reviewed and edited the content as needed and take full responsibility for the content of the publication. Built-in tools (such as GitHub Co-Pilot) have been used in the coding phase. 

\paragraph{\textbf{Data availability statement}} The data that support the findings of this study are available from the corresponding author upon reasonable request. They will sdeposited in a public repository after after the end of the research project. 

\paragraph{\textbf{Acknowledgements}} The authors thank Vivina Barutello, Susanna Terracini, Mattia G. Bergomi, Pietro Vertechi, Davide L. Ferrario. 

\paragraph{\textbf{CRediT authorship contribution statement}} \textbf{Diego Berti:} Writing – original draft, Visualization, Validation, Supervision, Software, Resources, Project administration, Methodology, Investigation, Funding acquisition, Formal analysis, Data curation, Conceptualization, Writing – review \& editing. \textbf{Gian Marco Canneori:} Writing – original draft, Visualization, Validation, Supervision, Software, Resources, Project administration, Methodology, Investigation, Funding acquisition, Formal analysis, Data curation, Conceptualization, Writing – review \& editing. \textbf{Roberto Cic\-carelli:} Writing – original draft, Visualization, Validation, Supervision, Software, Resources, Project administration, Methodology, Investigation, Funding acquisition, Formal analysis, Data curation, Conceptualization, Writing – review \& editing. \textbf{Irene De Blasi:} Writing – original draft, Visualization, Validation, Supervision, Software, Resources, Project administration, Methodology, Investigation, Funding acquisition, Formal analysis, Data curation, Conceptualization, Writing – review \& editing. \textbf{Margaux Introna:} Writing – original draft, Visualization, Validation, Supervision, Software, Resources, Project administration, Methodology, Investigation, Funding acquisition, Formal analysis, Data curation, Conceptualization, Writing – review \& editing. \textbf{Davide Polimeni:} Writing – original draft, Visualization, Validation, Supervision, Software, Resources, Project administration, Methodology, Investigation, Funding acquisition, Formal analysis, Data curation, Conceptualization, Writing – review \& editing.

\bibliography{references}
\bibliographystyle{plain}

\medskip

\noindent
Diego Berti \\
Dipartimento di Matematica ``Giuseppe Peano'', Universit\`a degli Studi di Torino\\
Via Carlo Alberto 10, 10123 Torino, Italy\\
\texttt{diego.berti@unito.it}

\medskip

\noindent
Gian Marco Canneori \\
Dipartimento di Matematica ``Giuseppe Peano'', Universit\`a degli Studi di Torino\\
Via Carlo Alberto 10, 10123 Torino, Italy\\
\texttt{gianmarco.canneori@unito.it}

\medskip

\noindent
Roberto Ciccarelli \\
Dipartimento di Matematica ``Giuseppe Peano'', Universit\`a degli Studi di Torino\\
Via Carlo Alberto 10, 10123 Torino, Italy\\
\texttt{roberto.ciccarelli@unito.it}

\medskip

\noindent
Irene De Blasi \\
Dipartimento di Matematica ``Giuseppe Peano'', Universit\`a degli Studi di Torino\\
Via Carlo Alberto 10, 10123 Torino, Italy\\
\texttt{irene.deblasi@unito.it}

\medskip

\noindent
Margaux Introna \\
Dipartimento di Matematica ``Giuseppe Peano'', Universit\`a degli Studi di Torino\\
Via Carlo Alberto 10, 10123 Torino, Italy\\
Department of Physics, University of Trento
Via Sommarive 14, 38123 Povo, Italy
\texttt{margaux.introna@unitn.it}

\medskip

\noindent
Davide Polimeni \\
Dipartimento di Matematica ``Giuseppe Peano'', Universit\`a degli Studi di Torino\\
Via Carlo Alberto 10, 10123 Torino, Italy\\
\texttt{davide.polimeni@unito.it}

\end{document}